\documentclass[letterpaper,oneside,10pt]{article}

\usepackage[USenglish]{babel}
\usepackage[T1]{fontenc}
\usepackage[ansinew]{inputenc}

\usepackage{lmodern} 

\usepackage{graphicx} 

\usepackage{float}
\usepackage{caption}
\usepackage{amsmath}
\usepackage{amsthm}
\usepackage{amsfonts}
\usepackage{amssymb}
\usepackage{mathtools}
\usepackage{mathdots}
\usepackage[hyphens]{url}

\let\Right\right 
\let\Left\left 
\makeatletter 
\def\right#1{\Right#1\@ifnextchar){\!\right}{}} 
\def\left#1{\Left#1\@ifnextchar({\!\left}{}} 
\makeatother

\newtheorem{theorem}{Theorem}
\newtheorem{lemma}{Lemma}

\begin{document}

\pagestyle{empty} 


\title{A matrix variation on Ramus's identity\\ for lacunary sums of binomial coefficients}
\author{John Blythe Dobson (j.dobson@uwinnipeg.ca)}

\maketitle


\pagestyle{plain} 


\begin{abstract}
\noindent
We study the well-known lacunary sums of binomial coefficients considered, most notably, by Christian Ramus, and their connection to a special kind of harmonic number associated with the first case of Fermat's Last Theorem. For one case of Ramus's famous identity we obtain a variation in which some of the parameters are replaced by square matrices of arbitrary dimension.

\end{abstract}

\section{Introduction}

\noindent
This paper examines the family of lacunary binomial sums introduced by Antoine Augustin Cournot in 1829 \cite{Cournot}, and considered to great effect by the Danish mathematician Christian Ramus (1806--1856) in 1834 \cite{Ramus}. These take the form

\begin{equation} \label{Definition}
T(N, r, m) := \sum_{\substack{j=0 \\ j \equiv r \bmod{N}}}^{m} \binom{m}{j} \quad (N > 1),
\end{equation}

\noindent
and have been the subject of an extensive literature, which will only be selectively reviewed here. Moreover, we do not attempt to treat the further multisections of these sequences into separate congruence classes of $m$ modulo $N$, which are usually encountered with $j, m \equiv 0 \pmod{N}$ and written in the form

\begin{displaymath}
S(N, 0, m) := \sum_{j=0}^{m} \binom{Nm}{Nj} = T(N, 0, N \cdot m);
\end{displaymath}

\noindent
nor do we attempt to treat the closely related \textit{jump sums} studied by Russell Jay Hendel (\cite{Hendel}, and OEIS A244608).

Ramus's remarkable identity, which reduces (\ref{Definition}) to a sequence in $N$, is

\begin{equation} \label{eq:Ramus}
T(N, r, m) = \frac{1}{N} \sum_{j=0}^{N-1} \left( 2 \cos\frac{j\pi}{N} \right)^m \cos\frac{j(m-2r)\pi}{N} \quad (m > 0).
\end{equation}

\noindent
It can also be written as

\begin{equation} \label{eq:RamusModernized}
T(N, r, m) = \frac{1}{N} \sum_{j=0}^{N-1} \omega^{(N-j)r} \left( 1 + \omega^j \right)^m \quad (m > 0),
\end{equation}

\noindent
where $\omega = e^{\frac{2\pi{}i}{N}}$ is a primitive $N$th root of unity. This formula is used by Hoggatt and Alexanderson \cite{HoggattAlexanderson} to derive closed-form expressions for $N = 2, 3, 4, 5, 6, 8$, and by Howard and Witt \cite{HowardWitt} for $N=10$. A natural counterpart to this is the alternating sum

\begin{equation} \label{eq:SunDefinitionAlternating}
\begin{aligned}
T^{\ast}(N, r, m) & := \sum_{\substack{j=0 \\ j \equiv r \bmod{N}}}^{m} (-1)^{(j-r)/N} \binom{m}{j} \\
                  & = \frac{1}{N} \sum_{\substack{j=1 \\ j \text{ odd} }}^{2N - 1} \omega^{(N-j)r} \left( 1 + \omega^j \right)^m \quad (m > 0) \\
                  & = 2 \cdot T(2N, r, m) - T(N, r, m).
\end{aligned}
\end{equation}

\noindent
This is evaluated by Howard and Witt \cite{HowardWitt} in the cases of $N = 2, 3, 4, 5$.

Our method, for reasons that will be explained below, is only applicable to the case where $r = 0$, and with $\omega$ as above,

\begin{equation} \label{eq:RamusSpecialCase}
T(N, 0, m) = \frac{1}{N} \sum_{j=0}^{N-1} \left( 1 + \omega^j \right)^m \quad (m > 0)
\end{equation}

\noindent
and

\begin{equation} \label{eq:RamusSpecialCaseAlternating}
\begin{split}
T^{\ast}(N, 0, m) & = \frac{1}{N} \sum_{\substack{j=1 \\ j \text{ odd} }}^{2N - 1} \left( 1 + \omega^j \right)^m \quad (m > 0) \\
                  & = 2 \cdot T(2N, 0, m) - T(N, 0, m).
\end{split}
\end{equation}

\noindent
Thus, our result is not a full generalization of Ramus's. However, it does generalize (\ref{eq:RamusSpecialCase}) in a different direction from (\ref{eq:Ramus}) or (\ref{eq:RamusModernized}), extending the solutions in simple radicals from the instances $N$ =  2, 3, 4, 5, 6, 8, 10 to all $N$ that are \textit{multiples} of 2, 3, 4, 5, 6, 8, or 10.

As a preliminary to developing this result, we note that the work of Konvalina and Liu (\cite{KonvalinaLiu}, p. 12) implies that $T(N, 0, m)$ (with $N$ even) and $T^{\ast}(N, 0, m)$ (with $N$ odd) satisfy the same simple recurrence relations of order $N-1$,

\begin{equation} \label{eq:CharacteristicPolynomialT}
\sum_{j=0}^{N-1} (-1)^j \binom{N}{j}a_{n-j} = 0.
\end{equation}

\noindent
Likewise, $T(N, 0, m)$ (with $N$ odd) and $T^{\ast}(N, 0, m)$ (with $N$ even) satisfy recurrence relations also of order $N-1$,

\begin{equation} \label{eq:CharacteristicPolynomialTStar}
\pm{}2a_{n-N} + \sum_{j=0}^{N-1} (-1)^j \binom{N}{j}a_{n-j} = 0.
\end{equation}

\noindent
A few examples of the coefficients in these recurrences are given in the Appendix, along with some values of the corresponding sequences.

\section{Matrix representations}

\noindent
Charles S. Kahane \cite{Kahane} has shown that the characteristic polynomial corresponding to the recurrence for $T(N, 0, m)$ (\ref{eq:CharacteristicPolynomialT}) is embodied in the simple circulant matrix of dimension $N$,

\begin{equation*}
C_N := \text{Circ}_N(1, 1, 0, 0, 0, \ldots) = 
\begin{pmatrix*}[r]
  1      & 1      & 0      & 0      & \cdots & 0      & 0 \\
  0      & 1      & 1      & 0      & \cdots & 0      & 0 \\
  0      & 0      & 1      & 1      & \ddots & 0      & 0 \\
  \vdots & \vdots & \vdots & \ddots & \ddots & \vdots & \vdots \\
  0      & 0      & 0      & \cdots & 1      & 1      & 0 \\
  0      & 0      & 0      & \cdots & 0      & 1      & 1 \\
  1      & 0      & 0      & \cdots & 0      & 0      & 1
\end{pmatrix*}.
\end{equation*}

\noindent
This matrix is just the sum of the $N$-dimensional identity matrix $I_N$ and the $N$-dimensional special matrix $U_N := \text{Circ}_N(0, 1, 0, 0, 0, \ldots)$, variously known as the \textit{basic circulant permutation matrix} or \textit{primary circulant permutation matrix} or \textit{forward shift permutation matrix}, or as we shall call it the \textit{unit circulant matrix}. It may be readily deduced from Kahane's result that the characteristic polynomial of $T^{\ast}(N, 0, m)$ is represented by the skew version of this matrix, $C_N^{\ast} := \text{Circ}_N^{\ast}(1, 1, 0, 0, 0, \ldots)$, which is the same as $C_N$ except that all the values below the diagonal --- in this case the lone 1 in the lower-left --- are negated. But these two matrices have the additional property, crucial for our purposes, that the first terms of the first rows of their powers exactly coincide with the terms of the respective sequences, a fact that will be proved below. Not being aware of a standard notation for such a relation, we henceforth use the symbol $\doteq$, hoping to avoid any possible confusion with ordinary equality. Some examples for small $N$ will illustrate the pattern:

\begin{displaymath}
T(2, 0, m) \doteq
\begin{pmatrix*}[r]
1 & 1 \\
1 & 1
\end{pmatrix*}^m
\qquad
T^{\ast}(2, 0, m) \doteq
\begin{pmatrix*}[r]
 1 & 1 \\
-1 & 1
\end{pmatrix*}^m
\end{displaymath}

\begin{displaymath}
T(3, 0, m) \doteq
\begin{pmatrix*}[r]
 1 & 1 & 0 \\
 0 & 1 & 1 \\
 1 & 0 & 1
\end{pmatrix*}^m
\qquad
T^{\ast}(3, 0, m) \doteq
\begin{pmatrix*}[r]
 1 & 1 & 0 \\
 0 & 1 & 1 \\
-1 & 0 & 1
\end{pmatrix*}^m.
\end{displaymath}

\noindent
In order for the powers of a matrix to coincide with the terms of a recurrence of degree $N$, a standard result of matrix theory states that they must agree in their characteristic polynomials and in their first $N$ terms (counting from term 0). Thus it only remains to demonstrate that the second of these conditions is met:

\begin{lemma} \label{Lemma_1}
For the $N$-dimensional circulant matrix $C_N$ representing $T(N, 0, m)$ and the $N$-dimensional skew-circulant matrix $C_N^{\ast}$ representing $T^{\ast}(N, 0, m)$, the first term of the first row of each of their first $N$ powers (starting from the power 0) is 1.

\begin{proof}
Let $U_N$ be the $N$-dimensional unit circulant matrix, $U_N^{\ast}$ be the skew version thereof, and $I_N$ be the $N$-dimensional identity matrix; then by definition

\begin{equation} \label{eq:PowerFormula}
C_N^k = \left\{ C_N + I_N \right\}^k = \sum_{j=0}^{k} \binom{k}{j} I^{n-j}_N \cdot U^j_N = \sum_{j=0}^{k} \binom{k}{j} U_N^j,
\end{equation}

\noindent
while

\begin{equation} \label{eq:PowerFormulaStar}
(C_N^{\ast})^k = \left\{ C_N^{\ast} + I_N \right\}^k = \sum_{j=0}^{k} \binom{k}{j} I^{n-j}_N \cdot (U_N^{\ast})^j = \sum_{j=0}^{k} \binom{k}{j} (U_N^{\ast})^j.
\end{equation}

\noindent
The values of $U^k_N$, starting at $k=0$, are

\begin{displaymath}
\begin{split}
& I_N,\\
& \text{Circ}_N(0, 1, 0, 0, 0, \ldots),\\
& \text{Circ}_N(0, 0, 1, 0, 0, \ldots),\\
& \text{Circ}_N(0, 0, 0, 1, 0, \ldots),
\end{split}
\end{displaymath}

\noindent 
et cetera, until for $k=N$ they return to $I_N$. The values of $(U_N^{\ast})^k$, starting at $k=0$, are the same as these, with no appearance of negated terms \textit{in the first row} (or indeed anywhere but in the lower-left triangle), until for $k=N$ they reach $-I_N$. Thus for $k<N$ the final sum in both (\ref{eq:PowerFormula}) and in (\ref{eq:PowerFormulaStar}) corresponds to the circulant matrix

\begin{displaymath}
\text{Circ}_N\left(\binom{k}{0}, \binom{k}{1}, \binom{k}{2}, \ldots, \binom{k}{k}[, \ldots]\right),
\end{displaymath}

\noindent
with leading term $\binom{k}{0} = 1$, and with the bracketted expression indicating trailing zeros when $k<N-1$. Thus the first $N$ values of $U^k_N$ or of $(U_N^{\ast})^k$ (counting in either case from $k=0$) supply the required $N$ consecutive first terms of first rows equal to 1. Hence the Lemma follows.

\end{proof}
\end{lemma}

We are now ready to state:

\begin{theorem} \label{Theorem_1}
For the $N$-dimensional circulant matrix $C_N$ and skew-circulant matrix $C_N^{\ast})$, the following relations hold:

\begin{equation} \label{CirculantFormula}
\begin{split}
       T(N, 0, m) & \doteq C_N^m; \\
T^{\ast}(N, 0, m) & \doteq (C_N^{\ast})^m.
\end{split}
\end{equation}

\begin{proof}
It follows from the work of Kahane \cite{Kahane} that the circulant matrix $C_N$ conforms in its characteristic polynomial with $T(N, 0, m)$, and that the skew-circulant matrix $C_N^{\ast}$ conforms in its characteristic polynomial with $T^{\ast}(N, 0, m)$. Lemma \ref{Lemma_1} establishes that the first $N$ values of the sequences agree in each case with the $N$ first terms of the first rows of the matrix powers, as required. Thus, the matrices represent the respective sequences in the sense defined in the beginning of this section.
\end{proof}

\end{theorem}

\section{A variation on Ramus's identity}

Henceforth, for the sake of simplicity, we confine most of the discussion to $T(N, 0, m)$, since $T^{\ast}(N, 0, m)$ can be so readily derived therefrom. As is well known, circulant and skew-circulant matrices whose dimension is a composite number divisible by $d$ can be partitioned into $d$ square blocks, and the matrix thus partitioned has the same powers as the original matrix when the operations are performed element-wise (\cite{Davis}, p. 18). For example, a matrix of dimension 6 can be partitioned in two different ways, with $d = 2, 3$:

\begin{equation*}
\begin{pmatrix*}[r]
   1 & 1 & 0 & \vdots & 0 & 0 & 0 \\
   0 & 1 & 1 & \vdots & 0 & 0 & 0 \\
   0 & 0 & 1 & \vdots & 1 & 0 & 0 \\
 \hdotsfor{7} \\
   0 & 0 & 0 & \vdots & 1 & 1 & 0 \\
   0 & 0 & 0 & \vdots & 0 & 1 & 1 \\
   1 & 0 & 0 & \vdots & 0 & 0 & 1
\end{pmatrix*}
\quad \text{or} \quad
\begin{pmatrix*}[r]
    1 & 1 & \vdots & 0 & 0 & \vdots & 0 & 0 \\
    0 & 1 & \vdots & 1 & 0 & \vdots & 0 & 0 \\
 \hdotsfor{8} \\
    0 & 0 & \vdots & 1 & 1 & \vdots & 0 & 0 \\
    0 & 0 & \vdots & 0 & 1 & \vdots & 1 & 0 \\
 \hdotsfor{8} \\
    0 & 0 & \vdots & 0 & 0 & \vdots & 1 & 1 \\
    1 & 0 & \vdots & 0 & 0 & \vdots & 0 & 1
\end{pmatrix*}.
\end{equation*}

\noindent
It is evident by inspection that a partition of such matrices into blocks is a compound matrix of dimension $d$ having the form

\begin{equation}  \label{eq:MatrixM}
M \medspace := \text{Circ}_d(a, b, 0, 0, 0, \ldots) = \medspace
\begin{pmatrix*}[r]
      a & b      & 0      & \hdots & 0 \\
      0 & a      & b      & \ddots & \vdots \\
      0 & 0      & a      & \ddots & 0 \\
 \vdots & \vdots & \vdots & \ddots & b \\
      b & 0      & 0      & \hdots & a
\end{pmatrix*},
\end{equation}

\noindent
where only $a$ and $b$ are non-zero, $a$ consisting of zeros except for a diagonal of 1s and a superdiagonal of 1s, and $b$ consisting of zeros except for a single $\pm{}1$ in the lower-left corner.

The general expression for the powers of a circulant compound matrix in which the submatrices are not expected to commute under multiplication is given in a important paper of 1969 by Lovass-Nagy and Powers \cite{LovassNagyPowers} which deserves to be better known. We require only the following very special case of an unnumbered result on p. 128 of their paper:

\begin{lemma} \label{Lemma_2}
For the compound matrix $M$ of dimension $d$ as in (\ref{eq:MatrixM}), with $a, b$ being square matrices and $\omega = e^{\frac{2\pi{}i}{d}}$ a primitive $d$th root of unity, we have

\begin{equation} \label{LovassNagyPowers}
M^m \doteq \frac{1}{d} \sum_{j=0}^{d - 1} \left( a + b \cdot \omega^j \right)^m \quad (m > 0).
\end{equation}

\end{lemma}

\noindent
As an immediate consequence, we have the following corollary (which we designate a theorem only for convenience of reference):

\begin{theorem} \label{Theorem_2}
If $N$ is divisible by $d$, the circulant matrix $C_N$ may be expressed as $\text{Circ}_d(a, b, 0, 0, 0, \ldots)$ as in (\ref{eq:MatrixM}), with $a, b$ being square matrices of dimension $N/d$ such that $a$ consists of zeros except for a diagonal of 1s and a superdiagonal of 1s, and $b$ consists of zeros except for a 1 in the lower-left corner; then for $\omega = e^{\frac{2\pi{}i}{d}}$ a primitive $d$th root of unity, 

\begin{equation} \label{eq:SplitGeneral}
T(N, 0, m) \doteq \frac{1}{d} \sum_{j=0}^{d - 1} \left( a + b \cdot \omega^j \right)^m \quad (m > 0).
\end{equation}

\end{theorem}

\noindent
The correspondence of this theorem with Ramus's identity (\ref{eq:RamusModernized}) is obvious in the extremal case $d = N$, which gives $a = 1$, $b = 1$. In the other extremal case, $d = 1$, it gives $a = C_N$, $b = 0$, and thus corresponds with Theorem \ref{Theorem_1}. For $d = 2$, the only other value of $d$ leading to a totally real matrix, it gives

\begin{equation} \label{eq:SplitTwo}
\begin{aligned}
T(N, 0, m) & \doteq \frac{(a + b)^m + (a - b)^m}{2}
           & = T(N/2, 0, m) + T^{\ast}(N/2, 0, m), 
\end{aligned}
\end{equation}

\noindent
providing an alternate derivation of (\ref{eq:SunFormulaSupplementEven}) below. Using this in conjunction with (\ref{eq:RamusSpecialCaseAlternating}), we further have for even $N$,

\begin{equation} \label{eq:SplitStar}
\begin{split}
T^{\ast}(N, 0, m) & \doteq \frac{(a + bi)^{2m} + (a - bi)^{2m}}{2} \\
                  & \doteq \Re \left( \left\{ a \pm bi \right\}^{2m} \right),
\end{split}
\end{equation}

\noindent
where in the second part the choice of sign is immaterial since $\Re (a + bi) = \Re (a - bi)$. For $N$ divisible by 4, the decomposition of $T(N, 0, m)$ with $d = 4$ is analogous:

\begin{equation} \label{eq:SplitFour}
T(N, 0, m) \doteq \frac{1}{4} \left\{ \left(a + b \right)^m + \left(a - b \right)^m + \left(a + bi \right)^m + \left(a - bi\right)^m \right\}.
\end{equation}

\noindent
Next in simplicity is the decomposition of $T(N, 0, m)$ for $N$ divisible by 3, with $d = 3$:

\begin{equation} \label{eq:SplitThree}
T(N, 0, m) \doteq \frac{1}{3} \left\{ \left(a + b \right)^m + \left(a + b\cdot\frac{-1 + i \sqrt{3}}{2} \right)^m + \left(a + b\cdot\frac{-1 - i \sqrt{3}}{2} \right)^m \right\}.
\end{equation}

\noindent
Setting

\begin{displaymath}
a = 
\begin{pmatrix*}[r]
 1 & 1 & 0 \\
 0 & 1 & 1 \\
 0 & 0 & 1
\end{pmatrix*},
\quad
b = 
\begin{pmatrix*}[r]
 0 & 0 & 0 \\
 0 & 0 & 0 \\
 1 & 0 & 0
\end{pmatrix*},
\end{displaymath}

\noindent
in (\ref{eq:SplitThree}) gives $(a + b)^m \doteq T(3, 0, m)$, so that

\begin{equation} \label{eq:ExpressionNine}
T(9, 0, m) \doteq \frac{1}{3} \cdot \frac{2^m + 2\cos(m\pi/3)}{3} + \frac{2}{3} \cdot \Re \left( \left\{a + b\cdot\frac{-1 \pm i \sqrt{3}}{2} \right\}^m \right),
\end{equation}

\noindent
where the choice of sign in the last term is immaterial. If we wish, we can replace the cosine with the closed-form evaluation (Schwatt \cite{Schwatt}, p. 179)

\begin{displaymath}
\cos(m\pi/3) = \frac{1}{4} (-1)^{\lfloor \frac{n+1}{3} \rfloor} \left( 3 + (-1)^{n + \lfloor \frac{n+1}{3} \rfloor} \right).
\end{displaymath}

\noindent
This result (\ref{eq:ExpressionNine}) may be compared with the evaluation of $N = 9$ given by Hendel \cite{Hendel}, and of the analogous $s^{\ast}(0, 9)$ by Sun (\cite{ZHSun1992}, pt. 2, Theorem 2.4). Admittedly, results like (\ref{eq:ExpressionNine}), while they formally entail smaller matrices than those given by the general formula (Theorem \ref{Theorem_1}), would likely have little value for computational purposes.

The fact that the first term in the right-hand side of (\ref{eq:SplitGeneral}) is $\frac{1}{d}(a+b)^m \doteq T(N/d, 0, m)$ proves that $d \cdot T(N, 0, m) - T(N/d, 0, m)$ represents an integer sequence. For example, the recurrence for $3 \cdot T(9, 0, m) - T(3, 0, m)$ is

\begin{equation} \label{NineRecurrence}
x_{n} - 6x_{n-1} + 15x_{n-2} - 19x_{n-3} + 12x_{n-4} - 3x_{n-5} + x_{n-6} = 0.
\end{equation}

\noindent
It can be shown that such sequences have order $N - N/d$, but as their orders are higher than those of the corresponding sequences studied by Hendel \cite{Hendel}, it does not seem worth the space to present a proof here.

The reason we cannot expect to extend Theorem \ref{Theorem_2} so as to provide a full generalization of Ramus's identity (\ref{eq:RamusModernized}) is that when $r \neq 0$, $T(N, r, m)$ evidently cannot be represented by a circulant matrix for $N > 2$. At least, it certainly has no such representations for $N = 3, 4$, as established by an exhaustive search. Fortunately, however, it is the case $r=0$ that has interesting parallels with the special harmonic numbers considered in the following section.

\section{The connection to some special harmonic numbers, and the first case of FLT}

\noindent
Before its complete proof by Andrew Wiles, a major result for the first case of Fermat's last theorem (FLT), that is, the assertion of the impossibility of $x^p + y^p = z^p$ in integers $x, y, z$, none of which is divisible by the exponent and where $p$ is a prime greater than 2, was the 1909 theorem of Wieferich that the exponent $p$ must satisfy the congruence

\begin{displaymath}
q_p(2) := \frac{2^{p-1} - 1}{p} \equiv 0 \pmod{p},
\end{displaymath}

\noindent
where $q_p(2)$ is known as the Fermat quotient of $p$, base 2. This celebrated result was to be generalized in many directions. One of these was the extension of the congruence to bases other than 2, the first such step being the proof of an analogous theorem for the base 3 by Mirimanoff in 1910. Another grew out of the recognition that some of these criteria could be framed in terms of certain special Harmonic numbers of the form

\begin{equation} \label{eq:Harmonic}
H_{\lfloor p/N \rfloor} := \sum_{j=1}^{\lfloor p/N \rfloor} \frac{1}{j}
\end{equation}

\noindent
for $N > 1$, and with $\lfloor \cdot \rfloor$ denoting the greatest-integer function. Furthermore, $H_{\lfloor p/N \rfloor}$ is the case $k=0$ of a sum studied by Lerch \cite{Lerch} and other writers,

\begin{equation} \label{eq:SkulaSum}
s(k, N) := \sum_{\substack{j=\lfloor\frac{kp}{N}\rfloor + 1\\ j \neq p}}^{\lfloor\frac{(k + 1)p}{N}\rfloor} \frac{1}{j},
\end{equation}

\noindent
where $\lfloor \cdot \rfloor$ denotes the greatest-integer function, and it is always assumed that $p$ is sufficiently large that $s(k, N)$ contains at least one element; the provision $j \neq p$ is necessary when $k + 1 = N$. In their landmark joint paper of 1995, Dilcher and Skula \cite{DilcherSkula} established the surprising result that the vanishing modulo $p$ of the sum in (\ref{eq:SkulaSum}) was a necessary criterion for the failure of the first case of FLT for the exponent $p$, for all $N$ from 2 to 46, and all corresponding values of $k$.

In order to delineate the connection with lacunary binomial sums, it is also convenient to define an alternating version of the sum just mentioned,

\begin{equation} \label{eq:SkulaSumAlternating}
s^{\ast}(k, N) := \sum_{\substack{j=\lfloor\frac{kp}{N}\rfloor + 1\\ j \neq p}}^{\lfloor\frac{(k + 1)p}{N}\rfloor} \frac{(-1)^j}{j}.
\end{equation}

\noindent
Now such sums as those in (\ref{eq:SkulaSum}) and (\ref{eq:SkulaSumAlternating}) are extremely intractable to calculate for large $p$. Thus when the twin brothers Zhi-Hong Sun and Zhi-Wei Sun showed, in separate papers published just over a decade apart, that they could be evaluated modulo $p$ in terms of simple functions of the lacunary binomial sums considered above, the result was not only of great theoretical interest but provided a powerful new computational technique. First, Z.-H. Sun (\cite{ZHSun1992}, pt.\ 1) showed that for a prime $p$,

\begin{equation} \label{eq:SunFormula}
\frac{N(1 - T(N, 0, p))}{p}  \equiv
\begin{dcases}
\medspace s(0, N) \pmod{p} & \text{if $N$ is even} \\
s^{\ast}(0, N) \equiv -s(1, 2N) \pmod{p} & \text{if $N$ is odd}.
\end{dcases}
\end{equation}

\noindent
Later, Z.-W. Sun (\cite{ZWSun2002}; see also \cite{ZWSun2008}) completed the analysis by showing that for a prime $p$,

\begin{equation} \label{eq:SunFormulaComplement}
\frac{N(1 - T^{\ast}(N, 0, p))}{p}  \equiv
\begin{dcases}
\medspace s(0, N) \pmod{p} & \text{if $N$ is odd} \\
s^{\ast}(0, N) \equiv -s(1, 2N) \pmod{p} & \text{if $N$ is even}.
\end{dcases}
\end{equation}

\noindent
We note in passing the advantageous decomposition

\begin{equation} \label{eq:SunFormulaSupplementEven}
s(0, 2N) \equiv s(0, N) + s^{\ast}(0, N) \pmod{p},
\end{equation}

\noindent
where by symmetry the right-hand side has the same value regardless of the parity of $N$, a value in agreement with (\ref{eq:SplitTwo}), (\ref{eq:SunFormula}), and (\ref{eq:SunFormulaComplement}).

The insights supplied by the equivalences proved by the Sun brothers would allow all of the historical results for $s(N, 0)$ to be obtained in a uniform way. However, in actuality these have seemingly only been employed in Z.\,H.\ Sun's own paper of 1993 (\cite{ZHSun1992}, pt.\ 2, Theorems 2.1 and 2.4) to evaluate $s(0, 16)$ and $s^{\ast}(0, 9) [\equiv s(1, 18) \bmod{p}$]. The culmination of developments obtained by other techniques is represented by the 2015 study of Al-Shaghay and Dilcher \cite{AlShaghayDilcher}, which evaluated $s(k, 24)$ for all values of $k$, thus subsuming all previous results for which $N$ is a divisor of 24.

From (\ref{eq:SunFormula}) and (\ref{eq:SunFormulaComplement}), it will be obvious that $T(N, 0, p) \equiv T^{\ast}(N, 0, p) \equiv 1 \bmod{p}$ for all $N, p$, and that the Dilcher-Skula criteria for $p$ to fail the first case of Fermat's Last Theorem translate for sums of Ramus's type as

\begin{equation} \label{eq:FTLCriterion}
T(N, 0, p) \equiv 1 \pmod{p^2}
\end{equation}

\noindent
for all $N$ between 2 and 23, and all even $N$ between 24 and 46.

Finally, it may be mentioned that Emma Lehmer's famous result (\cite{Lehmer1938}, p.\ 360)

\begin{equation} \label{eq:Lehmer}
\binom{p-1}{\lfloor jp/N \rfloor} \equiv (-1)^{\lfloor jp/N \rfloor} \left(1 - pH_{\lfloor jp/N \rfloor} \right) \pmod{p^2}
\end{equation}

\noindent
allows evaluations of these special binomial coefficients to be applied to these special harmonic numbers. The important 2011 study by Kuzumaki and Urbanowicz \cite{KuzumakiUrbanowicz} found expressions in terms of simple linear recurrences for $N = 24$ and every value of $j < N$, and thus, implicitly, for $s(k, 24)$ for all values of $k$.

A comparison of (\ref{eq:Lehmer}) with (\ref{eq:SunFormula}) and (\ref{eq:SunFormulaComplement}) reveals the existence of a threefold correspondence between lacunary binomial sums of Ramus's type, the type of special harmonic numbers under discussion, and special binomial coefficients considered by Lehmer. However, while the special harmonic numbers have a natural relationship to the other two objects, the relationship of those two to each other would involve a very convoluted expression.

\section{Appendix}

\noindent
A few examples of the recurrences for the sequences $T(N, 0, m)$ and $T^\ast(N, 0, m)$ are given in Tables \ref{Table_1} and \ref{Table_2} below. These are given in the form

\begin{displaymath}
c_N \cdot x_N + c_{N-1} \cdot x_{N-1} + \ldots + c_0 \cdot x_0 = 0.
\end{displaymath}

\begin{table} [H]
\begin{center}
\caption{Recurrent sequences for $T(N, 0, m)$ for small values of $N$}
\label{Table_1}
\begin{tabular}{ r | l }
       $N$ & coefficients $c$ in descending order    \\
\hline
 2 & $1, -2, 0$ \\
 3 & $1, -3, 3, -2$ \\
 4 & $1, -4, 6, -4, 0$ \\
 5 & $1, -5, 10, -10, 5, -2$ \\
 6 & $1, -6, 15, -20, 15, -6, 0$ \\
 7 & $1, -7, 21, -35, 35, -21, 7, -2$ \\
 8 & $1, -8, 28, -56, 70, -56, 28, -8, 0$ \\
 9 & $1, -9, 36, -84, 126, -126, 84, -36, 9, -2$ \\
10 & $1, -10, 45, -120, 210, -252, 210, -120, 45, -10, 0$ \\
\end{tabular}
\end{center}
\end{table}

\begin{table} [H]
\begin{center}
\caption{Recurrent sequences for $T^\ast(N, 0, m)$ for small values of $N$}
\label{Table_2}
\begin{tabular}{ r | l }
       $N$ & coefficients $c$ in descending order    \\
\hline
 2 & $1, -2, 2$ \\
 3 & $1, -3, 3, 0$ \\
 4 & $1, -4, 6, -4, 2$ \\
 5 & $1, -5, 10, -10, 5, 0$ \\
 6 & $1, -6, 15, -20, 15, -6, 2$ \\
 7 & $1, -7, 21, -35, 35, -21, 7, 0$ \\
 8 & $1, -8, 28, -56, 70, -56, 28, -8, 2$ \\
 9 & $1, -9, 36, -84, 126, -126, 84, -36, 9, 0$ \\
10 & $1, -10, 45, -120, 210, -252, 210, -120, 45, -10, 2$ \\
\end{tabular}
\end{center}
\end{table}

Some examples of the two sequences with $m$ starting at 0 follow, with OEIS references provided where applicable:

\begin{itemize}

\item
$T(2, 0, m)$: 1, 1, 2, 4, 8, 16, 32, 64, 128, 256, 512, 1024, 2048, 4096, 8192, 16384, 32768, 65536, 131072, 262144, $\ldots$ (OEIS A011782)

\item
$T^{\ast}(2, 0, m)$: 1, 1, 0, --2, --4, --4, 0, 8, 16, 16, 0, --32, --64, --64, 0, 128, 256, 256, 0, --512, $\ldots$ (OEIS A146559)

\item
$T(3, 0, m)$: 1, 1, 1, 2, 5, 11, 22, 43, 85, 170, 341, 683, 1366, 2731, 5461, 10922, 21845, 43691, 87382, 174763, $\dots$ (OEIS A024493)

\item
$T^{\ast}(3, 0, m)$: 1, 1, 1, 0, --3, --9, --18, --27, --27, 0, 81, 243, 486, 729, 729, 0, --2187, --6561, --13122, --19683, $\ldots$ (OEIS A057681)

\item
$T(4, 0, m)$: 1, 1, 1, 1, 2, 6, 16, 36, 72, 136, 256, 496, 992, 2016, 4096, 8256, 16512, 32896, 65536, 130816, $\ldots$ (OEIS A038503)

\item
$T^{\ast}(4, 0, m)$: 1, 1, 1, 1, 0, --4, --14, --34, --68, --116, --164, --164, 0, 560, 1912, 4616, 9232, 15760, 22288, 22288, $\ldots$ (similar to OEIS A099586 which however is not defined for $m=0$)

\item
$T(5, 0, m)$: 1, 1, 1, 1, 1, 2, 7, 22, 57, 127, 254, 474, 859, 1574, 3004, 6008, 12393, 25773, 53143, 107883, $\ldots$ (OEIS A139398)

\item
$T^{\ast}(5, 0, m)$: 1, 1, 1, 1, 1, 0, --5, --20, --55, --125, --250, --450, --725, --1000, --1000, 0, 3625, 13125, 34375, 76875, 153750, $\ldots$ (not presently in OEIS)

\end{itemize}

\end{document}